\documentclass{article}

\usepackage{indentfirst}
\usepackage{amsmath,amsfonts,amsthm,amssymb}
\usepackage{mathrsfs}
\usepackage{amscd}

\def\co{\colon\thinspace}
\DeclareMathAlphabet{\mathsfsl}{OT1}{cmss}{m}{sl}

\newtheorem{thm}{Theorem}
\newtheorem{lem}[thm]{Lemma}
\newtheorem{cor}[thm]{Corollary}

\newtheorem*{thm*}{Theorem}

\theoremstyle{definition}

\newtheorem{rem}[thm]{Remark}

\begin{document}

\title{Addendum to: ``Knots, sutures and excision"}

\author{{Yi NI}\\{\normalsize Department of Mathematics, Massachusetts Institute of Technology}\\
{\normalsize 77 Massachusetts Avenue, Cambridge, MA
02139-4307}\\{\small\it Emai\/l\/:\quad\rm yni@math.mit.edu}}

\date{}
\maketitle

\begin{abstract}We observe that the main theorem in \cite{KMsuture}
immediately implies its analogue for closed $3$--manifolds.
\end{abstract}

\begin{thm}\label{HMClosedFibered}
Suppose $Y$ is a closed irreducible $3$--manifold, $F\subset Y$ is
a closed connected surface of genus $g\ge2$. If
$HM_{\bullet}(Y|F)\cong\mathbb Z$, then $Y$ fibers over the circle
with $F$ as a fiber.
\end{thm}

The case that $g=1$ is already treated in \cite{KMBook}, following
the argument of Ghiggini \cite{Gh}.

\begin{thm}{\rm\cite[Theorem 42.7.1]{KMBook}}\label{HMgenus1}
Suppose $Y$ is a closed irreducible $3$--manifold, $F\subset Y$ is
a torus, $\eta$ is a $1$--cycle in $Y$ that intersects $F$ once.
If $HM_{\bullet}(Y|F,\Gamma_{\eta})\cong\mathcal R$, then $Y$
fibers over the circle with $F$ as a fiber.
\end{thm}

\begin{rem}
The statement of \cite[Theorem 42.7.1]{KMBook} uses a field
$\mathbb K_{\eta}$ of characteristic $2$, because the proof
involves the surgery exact sequence whose proof requires
characteristic $2$. Kronheimer pointed out that this part can be
replaced by the Excision Theorem \cite[Theorem 3.2]{KMsuture},
which allows us to use any characteristic.
\end{rem}

Let $M$ be the manifold obtained by cutting $Y$ open along $F$.
 The two boundary components of $M$ are denoted by $F_-,F_+$. $M$ can be viewed as a sutured
manifold with empty suture.

\begin{lem}
$M$ is a homology product, namely, $$H_*(M,F_-)\cong
H_*(M,F_+)\cong0.$$
\end{lem}
\begin{proof}
By \cite{MengTaubes,Turaev,Kcommunication}, Turaev's torsion
function $T(Y,\mathfrak s)$ is, up to a sign, equal to the Euler
characteristic of $HM_{\bullet}(Y,\mathfrak s)$ when $b_1(Y)\ge2$
and $\mathfrak s$ is a non-torsion Spin$^c$ structure. The
argument in \cite[Section~3]{NiClosedMfd} shows that $M$ is a
homology product if $b_1(Y)\ge2$.

If $b_1(Y)=1$, as suggested by Kronheimer, one can consider the
double of $M$ along $\partial M$, denoted by $Z$. Of course
$b_1(Z)\ge2$. Moreover, by \cite[Theorem~3.1]{KMsuture} we have
$HM_{\bullet}(Z|F_+)\cong\mathbb Z$. Let $M_2$ be the double of
$M$ along $F_-$, then $M_2$ is a homology product as in the last
paragraph. Now \cite[Lemma~4.2]{NiClosedMfd} implies that $M$ is
also a homology product.
\end{proof}

\begin{lem}
Suppose $\{F=F_1,F_2,\dots,F_n\}$ is a maximal collection of
mutually disjoint, nonparallel, genus $g$ closed surfaces in $Y$,
such that each surface is homologous to $F$. $M_1,M_2,\dots,M_k$
are the components of the manifold obtained by cutting $Y$ open
along these surfaces, $\partial M_k=F_k\cup F_{k+1}$. Let
$\mathcal E_k$ be the subgroup of $H_1(M_k)$ spanned by the first
homologies of the product annuli in $M_k$. Then $\mathcal
E_k=H_1(M_k)$ for each $k$.
\end{lem}
\begin{proof}
Since $M$ is a homology product, we can glue its two boundary
components together by a homeomorphism to obtain a new manifold
$Z$ such that $Z$ has the same homology as $F\times S^1$.  If
$\mathcal E_k\ne H_1(M_k)$ for some $k$, then as in
\cite[Section~4]{NiClosedMfd} we can construction two smooth taut
foliations $\mathscr F_1,\mathscr F_2$ of $Z$, such that
$F_k,F_{k+1}$ are compact leaves of $\mathscr F_1,\mathscr F_2$,
and
$$c_1(\mathscr F_1)\ne c_1(\mathscr F_2).$$
It then follows that $$\mathrm{rank}\;HM_{\bullet}(Z|F)>1$$ by
\cite[Corollary~41.4.2]{KMBook}. By \cite[Corollary~4.8]{KMsuture}
we have $HM_{\bullet}(Y|F)\cong HM_{\bullet}(Z|F)$, which is a
contradiction to the assumption that
$HM_{\bullet}(Y|F)\cong\mathbb Z$.
\end{proof}

\begin{cor}\label{CharSurj}
Let $(\Pi_k,\Psi_k)$ be the characteristic product pair (see
\cite[Definition~6]{NiCorr}) for $(M_k,\partial M_k)$, then the
map
$$i_*\co H_1(\Pi_k)\to H_1(M_k)$$ is surjective.
\end{cor}
\begin{proof}
See the proof of \cite[Corollary 7]{NiCorr}.
\end{proof}

\begin{proof}[Proof of Theorem \ref{HMClosedFibered}]
By Corollary~\ref{CharSurj}, each $\Pi_k$ contains a submanifold
$G_k\times I$, where $G_k$ is a genus $1$ surface with one
boundary component. Cutting $Y$ open along $F_k$'s and regluing by
suitable homeomorphisms, we can get a new manifold $Y'$ such that
the $G_k\times I$'s match together to form a submanifold $G\times
S^1\subset Y'$, where $G$ is a genus $1$ surface with one boundary
component. By \cite[Corollary~4.8]{KMsuture}, we have
$$HM_{\bullet}(Y'|F)\cong HM_{\bullet}(Y|F)\cong\mathbb Z.$$

Let $M'$ be the manifold obtained by cutting $Y'$ open along $F$,
then $M'$ is a homology product, and $M'$ contains a product
submanifold $G\times I$. Let $M''$ be the exterior of $G\times I$
in $M'$, and let $\gamma=(\partial G)\times I$. Then
$(M'',\gamma)$ is a sutured manifold which is a homology product.
By \cite[Definition~4.3]{KMsuture}, we have
$$SHM(M'',\gamma)= HM_{\bullet}(Y'|F)\cong\mathbb Z.$$
Now \cite[Theorem~6.1]{KMsuture} implies that $M''$ is a product,
thus $M'$ is also a product. So $Y'$ and hence $Y$ fiber over the
circle.
\end{proof}

\noindent{\bf Acknowledgements.}\quad We are very grateful to
Peter Kronheimer for lots of comments which improve this note
greatly. The author is partially supported by an AIM Five-Year
Fellowship and NSF grant number DMS-0805807. This note is written
when the author visited Peking University. The author wishes to
thank Shicheng Wang for his hospitality during the visit.


\begin{thebibliography}{H}

\bibitem{Gh}{ P Ghiggini}, {\it Knot Floer homology detects genus-one fibred
knots}, to appear in Amer. J. Math., available at
arXiv:math.GT/0603445

\bibitem{Kcommunication} P Kronheimer, private communication.

\bibitem{KMBook} P Kronheimer, T Mrowka, {\it Monopoles and
three-manifolds}, Cambridge University Press (2007).

\bibitem{KMsuture} P Kronheimer, T Mrowka, {\it Knots, sutures and
excision}, preprint, available at arXiv:0807.4891.

\bibitem{MengTaubes} G Meng, C Taubes, {\it $\underline{SW}=$ Milnor
torsion,} Math. Res. Lett. 3 (1996), 661--674.


\bibitem{NiCorr}{ Y Ni}, {\it Corrigendum to ``Knot Floer homology detects
fibred knots"}, preprint, available at arXiv:0808.0940.

\bibitem{NiClosedMfd} Y Ni, {\it Heegaard Floer homology and fibred
$3$--manifolds}, preprint, available at arXiv:0706.2032.

\bibitem{Turaev} V Turaev, {\it A combinatorial formulation for the Seiberg--Witten invariants of
$3$--manifolds}, Math. Res. Lett. 5 (1998), 583--598.

\end{thebibliography}
\end{document}